# Powers and division in the 'mathematical part' of Plato's *Theaetetus*

**Luc Brisson**  
Centre Jean Pépin  
CNRS-UMR 8230

**Salomon Ofman**  
CNRS-Institut mathématique de Jussieu-Paris Rive Gauche  
Histoire des Sciences mathématiques

*Abstract*.- In two articles ([Brisson-Ofman1, 2]), we have analyzed the so-called 'mathematical passage' of Plato's *Theaetetus*, the first dialogue of a trilogy including the *Sophist* and the *Statesman*. In the present article, we study an important point in more detail, the 'definition' of 'powers' ('δυνάμεις'). While in [Brisson-Ofman2], it was shown that the different steps to get the definition are mathematically and philosophically incorrect, it is explained why the definition itself is problematic. However, it is the first example, at least in the trilogy, of a definition by division. This point is generally ignored by modern commentators though, as we will try to show, it gives rise, in a mathematical context, to at least three fundamental questions: the meaning(s) of '*logos*', the connection between 'elements and compound' and, of course the question of the 'power(s)'. One of the main consequences of our works on Theaetetus' 'mathematical passage', including the present one, is to challenge the so-called 'main standard interpretation'. In particular, following [Ofman2014], we question the claim that Plato praises and glorifies both the mathematician Theodorus and the young Theaetetus. According to our analysis, such a claim, considered as self-evident, entails many errors. Conversely, our analysis of Theaetetus' mathematical mistakes highlights the main cause of some generally overlooked failures in the dialogue: the forgetting of the '*logos*', first in the 'mathematical part', then in the following discussion, and finally the failure of the four successive tries of its definition at the end of the dialogue. Namely, as we will show, the passage is closely connected with the problems studied at the end of the dialogue, but also to the two other parts of the trilogy through the method of 'definition by division'.

Finally, if our conclusions are different from the usual ones, it is probably because the passage is analyzed, maybe for the first time, simultaneously from the philosophical, historical and mathematical points of view. It had been considered usually either as an excursus by historians of philosophy (for instance [Burnyeat1978]), or as an isolated text separated from the rest of the dialogue by historians of mathematics (for instance [Knorr1975]), or lastly as a pretext to discuss some astute developments in modern mathematics by mathematicians (for instance [Kahane1985]).

*Résumé.- Dans deux articles en cours de rédaction* ([Brisson-Ofman1, 2])*, nous avons étudié ce qu'on appelle le 'passage mathématique' du* Théétète *de Platon, le premier dialogue d'une trilogie comprenant aussi le* Sophiste *et le* Politique*. Dans le présent article, nous considérons plus en détail un point important : la 'définition' des puissances (*'δυνάμεις'*), terme utilisés plusieurs fois avant d'être défini à la fin du passage. Dans [Brisson-Ofman2], nous avons montré que les différentes étapes pour obtenir cette définition sont incorrectes à la fois d'un point d'un vue mathématique et d'un point d'un vue philosophique. Nous expliquons ici pourquoi cette définition est elle-même problématiques. En dépit de ces erreurs, c'est la première tentative, au moins dans la trilogie, d'une définition par divisions, élément essentiel des deux autres dialogues. Ce point, souvent négligé par les commentateurs modernes, soulève, dans un cadre mathématique, certaines des plus importantes questions traitées dans la suite du dialogue : la signification du '*logos*', la relation 'éléments et composé', et la définition des 'puissances'. L'une des principales conséquences de nos études sur cette 'partie mathématique' est la mise en question de ce qu'on appelle 'l'interprétation* *standard principale'. En particulier, dans la perspective d'un précédent article ([Ofman2014]), nous montrons que la thèse standard selon laquelle le dialogue est l'éloge à la fois du mathématicien Théodore et du jeune Théétète, doit être révisée. Nous montrons inversement de quelle manière les erreurs mathématiques dans cette partie mettent en lumière la cause principale, généralement ignorée, des échecs successifs de Théétète : l'oubli du 'logos', tout d'abord dans le 'passage mathématique', puis dans la discussion qui suit, et finalement dans les quatre essais successifs de sa définition. En fait, le passage est étroitement lié aux problèmes étudiés à la fin du dialogue, mais son importance provient également de sa relation à la méthode générale de 'définition par division', outils principal dans les deux autres parties de la trilogie. Enfin, si nos conclusions diffèrent de celles usuelles, c'est qu'il est, et cela peut-être pour la première fois, analysé simultanément des points de vue philosophique, historique et mathématique. Il est en effet considéré généralement soit comme une digression sans suite par les historiens de la philosophie (par exemple [Burnyeat1978]), soit à étudier isolément du reste du dialogue par les historiens des mathématiques (ainsi [Knorr1975]) ou encore comme prétexte à discussions sur certains développement subtils en mathématiques modernes par les mathématiciens (ainsi [Kahane1985]).*

**Introduction.** Plato's three dialogues, the *Theaetetus*, the *Sophist* and the *Statesman* are explicitly connected at the beginning of the *Statesman*. [1] Many of their characters are the same: Socrates, the mathematician Theodorus from Cyrene, a Greek settlement in the north of actual Libya, Theaetetus a young Athenian, and a friend of him, *Socrates,* namesake of the first one. [2] He is silent here, but will take over Theaetetus in the *Statesman* (258b). The *Theaetetus* ends on Socrates' call for a next day meeting, while the *Sophist* opens on Theodorus reminding this announcement and introducing someone coming from Elea. [3] Socrates adds he would like to hear him about the sophist, the statesman and the philosopher (*Sophist*, 217a). Last, the *Statesman* begins by Theodorus' remark that, after the Sophist, the Stranger will search the definition of both the statesman and the philosopher, allowing Socrates to berate him as misapplied his mathematical knowledge, for he gives the 'equal worth (ἴσης ἀξίας)' to things which are 'beyond any proportion (ἀφεστᾶσιν ἢ κατὰ τὴν ἀναλογίαν)' (257a3-b7), a clear reference to the lesson on incommensurability in the 'mathematical part' of the *Theaetetus*.

In this dialogue the interlocutors try to find a definition of 'ἐπιστήμη' which may be understood either by 'science' or 'scientific knowledge', or even, as usually in English-language texts, simply by 'knowledge'. It is made clear from the very beginning that Socrates bears in mind the two first meanings, while Theaetetus the last one. Theaetetus' answer to the question 'what is science' consists in a list including some hand-crafted techniques. [4] Socrates criticized him not only for giving examples and not a definition, but because some of them, as skills, do not belong to the *definiendum*. In his answer, Socrates replaces them, without more explanation, by their associated theoretical sciences (146d-147a). Immediately afterwards begins the so-called 'mathematical passage'.

---

[1] Modern translations assume the three dialogues to be a trilogy (for instance [Klein1977], [Bernadete2008]). This chronology concerns only their literary presentation and entails nothing on the dating of their composition ([Klein1977], p. 3).

[2] To distinguish him from the older, its name will be written in italic.

[3] He is simply called the 'ξένος' in the dialogue. The usual translation is the 'Stanger', but it is also rendered by the 'Visitor' in some modern translations (for a discussion about this question, cf. [Sayre2006], note 1, p. 11). We follow here the usual translation.

[4] 'τε καὶ αἱ τῶν ἄλλων δημιουργῶν *τέχναί*' (146d1-2); cf. *Sophist* 221d, in particular 221d1-4, 231b, 267e, *contra* the *Statesman*, where they are used sometimes as synonymies, maybe because of the defects in the language (207d) or the different conceptions of the different interlocutors. For Mary Louise Gill these dialogues are some dialectical exercises for apprentice dialecticians ([Gill2016]). However, these questions would need a development beyond the purposes of the present article.

The analysis presented here concerns the method used to get the definition of 'powers' ('δυνάμεις'). It differs from the modern *doxa* or the so-called Plato's 'Main Standard Interpretation'. [5] As in our other studies, [6] we do not start from some supposed global interpretation of Plato's philosophy, but we take Plato's words seriously and stick to his words to try to understand the text the best we can, contrary to the modern '*hybris*' grounded of some claiming with Kant to 'better understand the text than [Plato] himself'. [7] The main difficulty for its understanding is the necessity, as in the whole passage, to study it simultaneously through three different fields: mathematics, philosophy and history.

1. **The 'powers'**

According to the 'Modern Standard Interpretation' or its shortcut MSI, [8] the 'mathematical passage' of the dialogue (147d-148b) is considered at best as paying some respect to Plato's friends, the two (actual or future) mathematicians Theodorus and Theaetetus, or even as being some distraction and rhetorical nicety with respect to its real philosophical content, the definition of science. In another article, [9] we tried to show that Plato's attitude is critical or at least ambivalent towards both Theodorus and Theaetetus, an important issue to understand the rest of the dialogue. As a matter of fact, nothing in the text supports the MSI claims. Namely, once again, it points to the exact opposite. Theodorus is presented as the kind of teachers filling the heads of his pupils but unable to develop their own power of thinking. [10]

Another common trend of modern interpretations is the claim that the exact meaning of the term 'powers'('δυνάμεις'), [11] both in mathematics and in Ancient Greek, is a crucial and obligatory step to understand the passage. [12] In [Brisson-Ofman2], we explained why this hunt for an elusive unique meaning is based on some misunderstandings. We want to show here it is an important term indeed, but mainly in order to appreciate Plato's critical approach of mathematics 'à la' Theodorus as well as the deep connection of the passage to the rest of the dialogue and to the two other parts of the trilogy.

---

[5] [Rowe2015], p. xiii. The literature about the trilogy is too large to be listed here. A good introduction for the *Theaetetus* and the *Sophist* will be found in [Rowe2015], p. xli-xlii and for the *Statesman* in [Brisson2003], p. 289-301.

[6] Cf. for instance [Brisson-Ofman1 and 2] and [Ofman2010 and 2014].

[7] 'when we compare the thoughts that an author expresses about a subject, in ordinary speech as well as in writings, it is not at all unusual to find that we understand him even better than he understood himself, since he may not have determined his concept sufficiently and hence sometimes spoke, or even thought, contrary to his own intention.' ([Guyer-Wood1998], Transcendental dialectic, I, 1, B370, p. 396).

[8] Though it is not part of it, since some scholars consider Fowler's book ([Fowler1999]) as some guide to understand the history of irrationality or even to understand the method used to get the incommensurability in *Theaetetus*' mathematical passage, let us make just two brief remarks. First according to the author himself, his book is highly speculative and has no historical basis (p. 18, 42, 369, 371, 373 ff.). As a matter of fact, it is a try to show that an important part of early Greek mathematics could have been obtained through the so-called 'anthyphairesis' (for a good analysis of it, cf. [Unguru2020]). Now concerning the *Theaetetus*, he writes himself that it was a too difficult task to try to connect this method to the dialogue (p. 381, note 24), and moreover, is not himself convinced that it was used by Theodorus or Theaetetus (p. 380). For a refutation of the 'anthyphairesis', as a supposed method used in *Thaetetus*, we refer to [Knorr1975], [Ofman2014] and [Brisson-Ofman1].

[9] [Brisson-Ofman2].

[10] These 'wise and inspired men' ('σοφοῖς τε καὶ θεσπεσίοις ἀνδράσι') to whom Socrates sends some of his aspiring disciples for a preliminary training (151b).

[11] Some scholars argue against this translation, but it is the usual one in ancient Greek.

[12] Arpád Szabó writes outright: 'my translation, as well as my interpretation of the text, is based primarily on the definition of a single word ['δύναμις']' ([Szabó1978], p. 36). Myles Burnyeat speaks about 'a vexed issue of terminology' and gives a large part of his analysis to this problem ([Burnyeat1978], p. 495-502; 'Among the most debated single terms of ancient Greek mathematics is the word *dynamis*' (J. Høyrup, "Dynamis…," *Historia Mathematica*, 17, 1990, p. 202). Cf. also, [Caveing1994], p. 172-173; [Knorr1975], p. 94; B. Vitrac, "Les formules de la « puissance »…", *Autour de la puissance chez Aristote*, M. Crubellier and al. (eds), Peeters, 2008, p.73-148.

*Theaetetus* (147d-148b) may be divided into two parts. The first concerns the account of Theodorus' mathematical lesson given for his young pupils (147d3-d6). [13] The second is about Theaetetus and his friend *Socrates* working together on the lesson, and trying to go further (147e5-148a4). [14] Let us give a brief summary of the first part.

Historians and commentators have fought for a long time on the term 'powers' ('δυνάμεις'). They can be essentially divided in two groups: for some, 'powers' would mean the 'sides' of squares, for others the 'squares' themselves. Contrary to the usual interpretation postulating a unique meaning, a thorough analysis of the passage leads to the conclusion that powers' ('δυνάμεις') has multiple meanings according to its different uses. [15] As the Stranger makes clear to Theaetetus in the *Sophist*, such kinds of ambiguities are the reason for the search of definitions: to get a mutual understanding on what the *definiendum* is (218b-c). [16]

- Its first instance, the subject of Theodorus' lesson ('Περὶ δυνάμεών τι ἡμῖν Θεόδωρος ὅδε ἔγραφε', 147d3), is vague: Theodorus drew (probably with a stick on the sand) some 'figures' representing the squares that involves both meanings, for the drawing of a figure shows both its sides and its area.
- The second ('ἐπειδὴ ἄπειροι τὸ πλῆθος αἱ δυνάμεις ἐφαίνοντο', 147d7-8) is about the infinite quantity of 'powers'. Once again, it is vague but its meaning is either the sides or the squares of areas non-perfect integers only, or eventually both, but in any case, different from the first one.
- The third ('ὅτῳ πάσας ταύτας προσαγορεύσομεν τὰς δυνάμεις', 147e1) introduces what Theaetetus and his friend *Socrates* learnt from Theodorus' lesson. It means beyond doubt the sides of the squares.
- The last one ('ὅσαι δὲ τὸν ἑτερομήκη, δυνάμεις' 148a7-8) concludes the boys' work. A definition of 'powers' is eventually given: they are the sides of the squares whose areas are of non-perfect integers.

It is not a sequence of more and more precise meanings ending in the definition; for instance, the third is not a consequence of the second. It is already a hint that the boys' work is not so methodical. [17]

Let us now consider how, using division, *Socrates* and Theaetetus come to the last, namely the only one, explicit definition of 'powers'.

## 2. The division

'Division' [18] plays a fundamental role in the *Sophist* as a tool of the dialectic method, as well as in the *Statesman* where it is conceptualized. [19] It is also used in the *Phaedrus* (266b) and in the *Philebus*. While there have been many studies of this method in these dialogues, not much is said about it in the *Theaetetus*, and more particularly its 'mathematical part', though it plays a fundamental role in the second part of the passage. Namely, it is the method used by Theaetetus and his friend *Socrates* to get a definition of the 'powers'.

---

[13] See [Brisson-Ofman1].

[14] See [Brisson-Ofman2]. Here we consider the particular point of the definition by divisions of the 'powers' that we did not detail in this previous article.

[15] This has been already argued, for instance in ([Allman1877], p. 271), but for different reasons, principally Plato's weaknesses in mathematics. Paul Tannery claimed the inconsistencies may result from a not well-defined mathematics vocabulary at Plato's time ([Tannery1876], I, p. 33) and also a series of mistakes by the copyists.

[16] Cf. the two kinds of 'love' in *Phaedrus*, one 'which is very justly reviled' ('ἐλοιδόρησεν μάλ᾽ ἐν δίκῃ') and another, a 'divine' one 'having the same name (ὁμώνυμον') as the first' (266a); or the two kinds of 'sophistry' in the *Sophist* (231b).

[17] The possibility that it is not so much the boys who will be judged but their master is left open by Socrates' exclamation at the end of the passage: 'I think Theodorus will not be found liable to an action for false witness.' (148b).

First, let us consider Plato's text itself. [20]

> THEAETETUS
> 
> (1) We divided the integer [21] in its totality in two parts. [22] One which has the power to be the product of an equal times an equal [23] [**148a**] we likened to the square as a figure, and we called it square and equilateral. [24]
> 
> SOCRATES
> 
> Good, so far.
> 
> THEAETETUS
> 
> (2) Then, the one in between, such as three or five or any one without the power to be a product of an equal times an equal, but is the product of a greater times a lesser, or a lesser times a greater, [25] that is always encompassed by a greater side and a smaller one, we likened it this time to a rectangular figure, [26] so that we called it a rectangular integer. [27]
> 
> SOCRATES
> 
> That's excellent. But how did you go on?
> 
> THEAETETUS
> 
> (3) All the lines squaring an equilateral and plane integer, we defined as 'length', while we defined as 'powers' the ones squaring the rectangular, because, although not commensurable as length with the formers, [**148b**] they are commensurable as areas that they have the power to produce. And it is the same in the case of solids. [28]

## 2.1. The question at stake

The boys considered the 'integer in its totality' ('ἀριθμὸν πάντα') and they 'divided it into two' ('δίχα διελάβομεν'). What does it mean?

---

[18] Namely, what is usually called 'division' in the literature means rather 'dichotomy' i.e. the particular division in two parts. It is connected to an important tool in arithmetic, the sequence of successive divisions by 2 of any integer (see [Ofman2010], p. 104-112).

[19] The term '*methodos*' appears twice in the *Sophist* and Aristotle refers to it as the 'method of division' ('ἡ διὰ τῶν διαιρέσεων ὁδός', *Post. An.* II, 5, 91b12).

[20] All the translations of Theaetetus' 'mathematical part' are from [Brisson-Ofman1,2].

[21] To avoid the confusion attached to the modern very large notion of 'number', we will always translate 'ἀριθμός' by 'integer'.

[22] Cf. *infra*, §2.2.

[23] The same sentence 'is the product of an equal times an equal' ('ἴσον ἰσάκις γίγνεσθαι') is found in Euclid's *Elements* (for instance, book VII, def. 18). It is a shortcut for saying that the integer is the product of the form $n \times n$, so that its result is a perfect square.

[24] 'Τὸν ἀριθμὸν πάντα δίχα διελάβομεν· τὸν μὲν δυνάμενον ἴσον ἰσάκις γίγνεσθαι τῷ τετραγώνῳ τὸ σχῆμα ἀπεικάσαντες τετράγωνόν τε καὶ ἰσόπλευρον προσείπομεν.'

[25] A shortcut (cf. also note 23, *supra*) to mean that the integer can only be written as a product of two different integers i.e. it is any integer except the perfect squares. Nevertheless, the terms 'equal' and 'unequal' are not names of objects but of relations (*x is* equal/unequal to *y*), so that it is a hint of the utmost importance of relations or ratios ('*logoi*') in mathematics.

[26] A rectangle in a strict sense that is with *unequal* sides i.e. *not* a square.

[27] 'Τὸν τοίνυν μεταξὺ τούτου, ὧν καὶ τὰ τρία καὶ τὰ πέντε καὶ πᾶς ὃς ἀδύνατος ἴσος ἰσάκις γενέσθαι, ἀλλ' ἢ πλείων ἐλαττονάκις ἢ ἐλάττων πλεονάκις γίγνεται, μείζων δὲ καὶ ἐλάττων ἀεὶ πλευρὰ αὐτὸν περιλαμβάνει, τῷ προμήκει αὖ σχήματι ἀπεικάσαντες προμήκη ἀριθμὸν ἐκαλέσαμεν.'

[28] "Ὅσαι μὲν γραμμαὶ τὸν ἰσόπλευρον καὶ ἐπίπεδον ἀριθμὸν τετραγωνίζουσι, μῆκος ὡρισάμεθα, ὅσαι δὲ τὸν ἑτερομήκη, δυνάμεις, ὡς μήκει μὲν οὐ συμμέτρους ἐκείναις, τοῖς δ'ἐπιπέδοις ἃ δύνανται. καὶ περὶ τὰ στερεὰ ἄλλο τοιοῦτον.'

The sentence (1) has to be connected to the long analysis at the end of the dialogue on the relation between the 'whole' and the 'parts' or the 'compound' and its 'elements' (202e-207c). Socrates emphasizes there the importance of the question, requesting Theaetetus' full attention. [29] First he asks him if he agrees that the 'consummate arithmetician knows *all the integers* (τι πάντας ἀριθμοὺς)', adding 'because he has the science of *all integers* in his soul (πάντων γὰρ ἀριθμῶν εἰσιν αὐτῷ ἐν τῇ ψυχῇ ἐπιστῆμαι)' (198b). He gets a strongly positive answer from Theaetetus.

Though generally unnoticed, such a statement is stunning, especially coming from someone considered as a brilliant student in mathematics. Since the integers are infinite in number, how would it be possible to *know* all the integers, *a fortiori* to *have* them in soul? [30]

### 2.2. The problems of the division

Then the discussion switches to the differences between 'the all' ('τὸ πᾶν'), 'the totality' ('τὰ πάντα') and 'the whole' ('τὸ ὅλον') (204a-205a), Socrates asking ('πότερον') Theaetetus if the 'all' is identical ('αὐτόν') or 'different' ('ἕτερον') of the 'whole' (204a-e). This time, a much less assertive Theaetetus chooses tentatively the second, changing his mind once again at the end of a long questioning by Socrates. He finally agrees that all these three terms ('the all', 'the totality', 'the whole') are the same. [31]

This development starts from the hypothesis that any 'perfectly scientific knowledge' ('τελείως (…) ἐπιστήμην') is a compound of elements that are 'irrational and unknowable' ('ἄλογα καὶ ἄγνωστα')' and can only 'be perceived through the senses' ('αἰσθητά') (202d). Moreover, it entails forcefully that the 'integer in its totality' is, for Theaetetus, nothing more than all its elements i.e. all the integers. This identification corresponds to Socrates' example, at the end of the dialogue, of the wagon which is certainly more than 'hundred pieces of wood' (207a-208c). What is lost is the layout between them i.e. their 'relations' ('*logoi*'). [32]

However, is it not true that the 'integer in its totality' is just the union of all the integers (a modern would speak of the set of all the integers)? A positive answer would forget the ordering, so that the integers would be arbitrarily given. Thus, for instance, Theodorus' sequence of odd integers from 3, 5 to 17 could as well contain 10001 (a *myriad* plus one) than 7. [33] Such an oversight will impede the boys' next

---

[29] 'τῷ δὲ δὴ ἐντεῦθεν ἤδη πρόσσχες τὸν νοῦν' (198b).

[30] Let us emphasize that it does not mean that it is impossible to know *some* properties on the *whole* integers. For instance, we know that all the integers are either odd or even. But it is evidently impossible to know *all* the integers: for instance, it is already impossible to know (at least till now) all the divisors of any integer chosen at random, and it is certainly impossible to have them in a human mind, not even in the huge memory of modern computers.

[31] The first definitive answer of Theaetetus to Socrates is that the combination is the 'complete totality' ('τὰ ἅπαντα ἔμοιγε δοκοῦμεν.', 203c). Only when Socrates shows him that this answer leads to, as the boy says, 'something monstrous and absurd' ('ἀλλὰ δεινὸν καὶ ἄλογον'), he accepts to change his mind once again (203e). Then, he reluctantly agrees to consider the other branch of the alternative 'because perhaps that will be better than the other way.' ('καὶ τάχα γ᾽ ἂν μᾶλλον οὕτως ἢ 'κείνως ἔχοι.', 203e), to conclude : 'I am not sure; but you tell me to answer boldly, so I take the risk and say that they are different.' ('ἔχω μὲν οὐδὲν σαφές, ὅτι δὲ κελεύεις προθύμως ἀποκρίνασθαι, παρακινδυνεύων λέγω ὅτι ἕτερον.', 204b). For a different analysis, cf. for instance [Harte2002], p. 40-41.

[32] As a matter of fact, to get a scientific knowledge of what a wagon is, the enumeration of its parts is not enough: the '*logos*' (207c) is needed. And the end of the dialogue is entirely about the meaning of this term. Conversely, Theodorus proudly says, previously in the dialogue, that 'rather too soon', he 'bent away' from the '*logoi*' ('ἡμεῖς δέ πως θᾶττον ἐκ τῶν ψιλῶν λόγων πρὸς τὴν γεωμετρίαν ἀπενεύσαμεν.', 165a).

[33] The examples used later are mostly about letters and syllables (as for instance 'S', 'O' and 'SO', 203a-e), the same words in Greek than respectively for 'elements' and 'compounds'. Then the necessity of an ordering (or

definition, but it entails also an ambiguity in the first sentence, the division of 'the integer in its totality'. As the compound ('the integer') and its elements ('the integers') are homonymous, the sentence may mean:
- o   The set of all integers is divided into two subsets.
- o   Each and every integer is divided by two.

Nevertheless, the following sentences will make it clear that Theaetetus means the first: another hint about the non-rigorous character of the boys' reflections.

This is another hint about the boys' awkwardness: the very beginning of their 'definition' is ambiguous i.e. not well-defined.

## 3.  Arithmetic and geometry

'Power' ('*dynamis*') is used here by the boys to relate arithmetic and geometry. It plays also a central role in Plato's doctrine of 'being', namely to link intellection and sense perception in the *Sophist* (246c-247e).

### 3.1. A first correspondence

The boys make an important correspondence, between arithmetic and geometry. Leaded by the figurative names given to the class or to their elements, the 'square and equilateral' integer, they associate to any of the perfect square integers a linear figure, the square of area this integer However, while in Theodorus' lesson, all the magnitudes, lines or surfaces, are given as a certain number of feet, here no unit is specified. Only abstract integers are taken into consideration by the boys, so that their idea needs to be that any integer $n$ is represented by a finite line, namely $n$ times an arbitrary (segment of) line $u$ considered as the length unit. [34] Then, to $n$ 'having the power to be a product of an equal times an equal' i.e. $n$ is a perfect square or in a modern representation $n = m^2$ for some integer $m$, they associate the square of side $m$. [35]

### 3.2. The figures

The definition by divisions into two parts or dichotomies is largely used in the *Sophist* and the *Statesman*. It is the privileged way to get a definition (i.e. what a thing is, *Sophist* 285d-e). The Stranger in the *Statesman* claims that 'by far our first and most important object should be to exalt the method itself of ability to divide by classes'. [36]

---

relation between the letters) is evident, since 'SO' is certainly different from 'OS'. It is different for the integers: the Ancient Greeks did not use a notation by position (cf. for instance in Plato's *Laws*, the integer 5040 is given as 'τεττᾰράκοντα καὶ πεντακισχιλίων' (738a) *contra* 'πεντακισχιλίων καὶ τεττᾰράκοντα' (771c)).

[34] As it is usually done in Euclid's *Elements*.

[35] More explicitly the side is of length $m$ times $u$. However, once the unit $u$ is fixed, it can be forgotten and we can say the length of the side is m. This is the common representation of integers (or more generally magnitudes) in Euclid's *Elements*. Socrate-Theaetetus' idea could be represented graphically in the figure below:

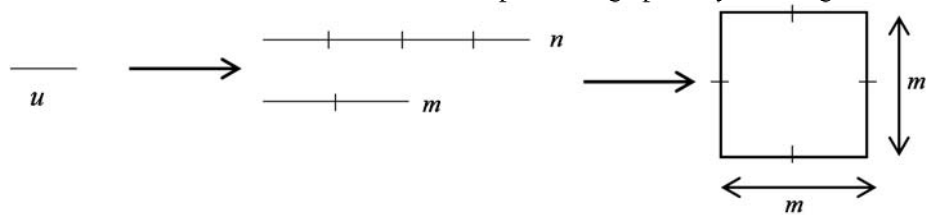

[36] 'πολὺ δὲ μάλιστα καὶ πρῶτον τὴν μέθοδον αὐτὴν τιμᾶν τοῦ κατ' εἴδη δυνατὸν εἶναι διαιρεῖν' (286d). The term 'εἶδος' is not easy to translate.

Thus, its use by the boys [37] seems to show that their method is correct. Let us consider it more in details.

### i) The squares

As the boys acknowledged previously (147d7, cf. *supra*, §1), the objects of their study are 'infinite in quantity', thus it is not possible to consider the problem 'one case in turn' as Theodorus did in his lesson (147d5). Hence, their goal here is to treat simultaneously an infinite quantity of cases.

As expected, the division gives two classes. The first is the 'square and equilateral' integer, defined as the one which has 'the power to be the product of an equal times an equal'. Theaetetus presents them rightly as 'having the power ('δυνάμενον') to be', rather simply 'to be', 'the product of an equal times an equal', for they are also always the product of two unequal integers (1 and the integer itself for example). For example, 4 is equal to the product of 1 and 4 as well as to the product of 2 and 2.

### ii) The Rectangles [38]

The definition of the second class, the one of 'rectangular integers' in sentence (2), raises some additional problems. Within 'the *integer* in its totality', the boys 'separate' ('διαλαμβάνειν', 147e) the integer which is a perfect square, [39] from the one that is not. This is problematic since none of the two rules laid down by the Stranger of the *Statesman* for a definitional division is respected.

The first is to get symmetrical parts. For instance, one cannot divide the 'mankind' by 'separating' ('ἀφαιροῦντες') Greeks and Barbarians, or the 'integer' by 'cutting off' ('ἀποτεμνόμενος') ten thousand from all the other integers (262d). A good division would be 'to divide 'the integer into odd and even', [40] or 'the mankind into male and female'. [41] The mistake in the definitional division is to 'take off one small part on its own, leaving many large ones behind, and without reference to real classes'. [42] Indeed, 'it's not safe to make thin cuts', [43] 'it's safer to go along cutting through the middle of things '. [44] Now, let us consider the boys' cutting-up (another sense of 'διαλαμβάνειν') the 'integer in its totality'. It is certainly neither symmetric, nor a division in half. As a matter of fact, between *two* successive perfect square integers there are *many* non-squares (for instance between $10\ 000 = (100)^2$ and $10201 = (101)^2$ there are 200 different non-squares). [45]

The second and most important rule is the unity of the resulting parts. They need to form a 'genus' ('γένος', 262d), i.e. to have some common characters. While Barbarians are precisely defined as non-Greeks, and all the integers minus ten thousand as 'non-ten thousand', the odd integers and the even integers have each a real unity, allowing a real definition, and it is the same with the male humans and the female humans. The only common character of all 'rectangular' integers is: *not* being 'square' integers. In other words, the division cut off the 'integer in its totality' into square and non-square, as in ten thousand and non-ten thousand, or the mankind off into Greek and non-Greek (Barbarian).

---

[37] However, the words are not the same in both cases: in the *Sophist*, it is 'διαιρεῖν', the same verb used in Euclid's *Element*, to define for instance the odd and the even integers, while Theaetetus uses the verb 'διαλαμβάνειν'.
[38] A rectangle here means a figure with *unequal* sides i.e. excluding the squares (cf. *supra*, note 26).
[39] 'τετράγωνόν τε καὶ ἰσόπλευρον' (148a).
[40] 'τὸν μὲν ἀριθμὸν ἀρτίῳ καὶ περιττῷ' (262e).
[41] 'τῶν ἀνθρώπων γένος ἄρρενι καὶ θήλει' (ib.).
[42] 'μὴ σμικρὸν μόριον ἓν πρὸς μεγάλα' (262b).
[43] 'ἀλλὰ γάρ, ὦ φίλε, λεπτουργεῖν οὐκ ἀσφαλές' (ib.).
[44] 'διὰ μέσων δὲ ἀσφαλέστερον ἰέναι τέμνοντας' (ib.).
[45] Though for modern elementary logics the asymmetry is meaningless because both sets are infinite, see *infra*, **Appendix**.

Another generally unnoticed problem is that the boys give two definitions of 'rectangular' integers. The first one is purely negative: they are not 'square' integers. The second is not really a definition but appears as some obviousness: they are the product of two unequal integers. Theaetetus emphasizes this point, adding they are products of 'a greater and a smaller'. [46]

The first part is obvious. Since they are not 'square' integers, they have not, by definition, the power to be a product of an 'equal times an equal' (i.e. to be a perfect square). But why are they a product of unequal integers? Namely, according to the popular *mantra* claiming that 1 was absolutely not an integer for the Ancient Greeks, this would be flatly false, for no prime integer could be considered as the product of two integers! [47] Thus, clearly, at least here for Plato, 1 **is** an integer, and any integer $n$ is indeed such a product i.e. $n = n \times 1$. [48]

Last, they repeat the correspondence done in point 3.1, *supra*. Once an arbitrarily (segment of) line $u$ set as the unit of length, to any 'rectangular integer' $n$ equal to the product of $p$ and $q$, they associate a rectangle of sides $p$ and $q$. [49]

### 3.3. True results, incorrect reasoning

Using the graphical names given to the two classes (or their elements) defined previously (the 'square and equilateral' and the 'rectangular'), they associate a figure to any integer ('square' or 'rectangular'). [50] However, while the figure associated to each element of the first class (the square integer) is well-defined (cf. *supra*, §3.1), the situation is different for the second class. Generally there are several possible rectangles associated to the same rectangular integer. For instance, since the products of 1 and 6, as well as of 2 and 3 gives 6. Thus, for any unit $u$ fixed, both the rectangles of sides (*1,6*) and (*2,3*) may be associated to the rectangular integer 6. The boys do not consider this problem. However, in sentence (3), they make an important connection between arithmetic and geometry.

i) For an integer $n$ which 'has the power to be the product of an equal times an equal', they associate the square of side $m$ (i.e. for $n = m^2$, the square of side $m$ where $m$ is an integer). They did not say it explicitly, but once again we have to suppose an arbitrary unit of length $u$ is fixed, and to the integer $n = m^2$ is associated the square of side $m$.

ii) For a 'rectangular integer' $n = p \times q$, they associate first a rectangle of sides $p$ and $q$, and then a square of same area as this rectangle. [51] As we saw, there is a problem for the choice of the rectangle. Nevertheless, the boys did not see it. Since they clearly know that there are in general many ways to get the same product, [52] so they thought, in their soul, they have a way to associate to each such integer a unique rectangle. Namely, in the same way as Theodorus drew the powers in his lesson, to each 'rectangular integer' *n*, they thought to associate the rectangle of sides (*1,n*) and then the square of area *n*. Since there is only one such rectangle, the problem of the choice between different rectangles disappears. This follows closely Theodorus' construction. [53]

---

[46] Cf. *supra*, notes 23 and 25.

[47] For instance, there would be no integers such that 3 would be equal to their product.

[48] Cf. also *infra*, §3.3. Double definitions are also sometimes found in Euclid's *Elements*, for instance in the definition VII.7, often blamed by historians of mathematics, of the odd integer which either cannot be divided in half or differs from an even integer by a unit.

[49] For the omission of u, see *supra*, note 35.

[50] Theodorus' lesson is based on some drawings, and the boys tried naturally to extend this method to any integer (cf. [Brisson-Ofman1]).

[51] Cf. supra, note 35.

[52] The case of 6 is considered in details in 204b10-c3.

Nevertheless, the boys did not see the problem of a possible ambiguity in the construction, that the reader has to fix by himself. It is the first of several cases where, while the results are true, the boys' reasoning is not correct. [54] This example shows that, though truth is certainly a property of science, it does not characterize it, as Socrates makes clear in his criticism of Theaetetus' proposal of his last definition of science as 'true opinion'. [55]

## 4. An 'awkward' definition

It is remarkable that the definitions of both 'length' and 'powers', given by *Socrates*-Theaetetus, disappeared from the field of mathematics, while the problem of the commensurability/incommensurability was still central for a very long period, at least till Euclid's time (cf. book X of the *Elements*). This is generally a hint of a problem and indeed, as we will see, the definition of the formers is incorrect, according to both Plato and the mathematicians of his time.

### 4.1. A too 'generous' division

First, the division into 'square' and 'rectangular' integers is simply useless. As a matter of fact, the boys would get the same result at once, if for any integer *n* ('square' or 'rectangular'):
- they had associated the rectangle of sides (1,*n*)
- and to this rectangle they had associated the square of same area as the rectangle. [56]

Once again, Theaetetus is too 'generous', [57] giving a too long entangled exposition. [58] But there is more to come.

### 4.2. A problematic reasoning

At the end of the passage, the boys define the 'powers' and the 'length' ('μῆκος'). According to the boys' definition, the latter is the side of the squares associated to the

---

[53] In [Brisson-Ofman1] we showed the simplest way for Theodorus to construct the 'powers' associated to 3 feet, 5 feet, etc. till 17 feet, was a follow: first to draw a rectangle of sides (1foot, *n* feet) (where *n* is equal to 3, 5, …, 17), then to use an immediate corollary of Pythagoras' theorem: in a right triangle BDO of hypotenuse OB and height DH, the rectangle on OH and HB is equal to the square on DH (in modern notations, OH×HB = DH$^2$).

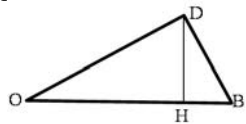

In particular, if we fix HB equal to the unit and OH to *n* times this unit, then DH is the side of the square equals (in area) to the rectangle (1 foot, *n* feet). For more details, we refer to [Brisson-Ofman1]. The boys generalized this method to any rectangular integer.

[54] For a much deeper mathematical problem about the same question, we refer to the **Appendix** in [Brisson-Ofman2].

[55] 'true judgment may well be knowledge. So let that be my answer' ('κινδυνεύει δὲ ἡ ἀληθὴς δόξα ἐπιστήμη εἶναι, καί μοι τοῦτο ἀποκεκρίσθω', 187b5-6).

[56] Cf. *supra*, note 53.

[57] He is 'φιλόδωρος' (146d) because, instead of giving 'a short (…) answer, [he went] an interminable distance round' (147c).

[58] Cf. previous paragraph.

'square' integers, in other words, the sides which are integers. [59] As the boys redefine 'length' in an unusual way, they also redefine 'powers': the sides of the squares of areas non-perfect square integers. In modern words, an element of the class 'length' is a length equal to an integer (for instance 3), and a 'power' is a length equal to a square root of any non-square integer (for instance √3). Since the term is here redefined, its previous meaning must have changed, supplementary evidence that its meaning cannot be unique in the passage.

Moreover, the redefinition of 'length' is surprising because Theodorus' lesson was precisely about lengths 'incommensurable to the unit' (147d5). However, the boys decided to call these *lengths* 'powers'. It seems senseless as long it is not linked to some parts of Euclid's *Elements*, in particular the definitions in book X. [60] In the *Elements*, the word 'powers' ('δυνάμεις') is not used to qualify some lines as in Theaetetus' account but as a property connected to a dative form ('δυνάμει', cf. for instance def. 2 and 3, as well as for the verb 'δύνασθαι' under the form 'δυνάμεναι' ('being able' to do/be something) as in definition 4. Definition 3 opposes lines 'incommensurable in length only' [61] to lines 'incommensurable in power'. [62] Moreover, this definition contains a result: the existence of an infinite quantity of lines 'incommensurable in length only'. [63] Thus, the boys changed the usual meaning of 'in length' and 'in power' for 'length' and 'powers'. In other words, the lines commensurable to the unit are defined as elements of the class 'length', while the incommensurable lines are defined as 'powers'. Plato shows, once again (cf. *supra*, note 25), the boys replacing a relation ('*logos*') between objects by the objects themselves.

Both definitions result from the division of the class of the 'integer in its totality'. But to which class do the elements of both classes ('length' and 'powers') belong together? Though Theaetetus did not name the elements of the 'length', we know they are, as well as the 'powers', sides of squares, hence lines. Thus, what we get finally are lines that are either an integer (i.e. the elements of the class 'length'), or that are not such, but 'have the power to produce' squares that are a whole number (i.e. an integer) of times the corresponding square-unit (i.e. the 'powers'). [64]

It is not easy to define such a mixed class because some of its elements come from arithmetic, others from geometry. In Plato's terms, something is missing: a common unity to form a '*genos*'.

---

[59] In modern notation, a 'square' integer is an integer $n$ of the form $m^2$, where $m$ is some integer. Thus the side of the square associated to $n$ is $\sqrt{m^2} = m$. Conversely, any integer $m$ is the side of the square associated to the perfect square integer $n = m^2$. Thus, once set a unit of length $u$, the 'length' according to the boys' meaning is formed of all the lines of length an integer (with respect to $u$).

[60] Thus these definitions came from much earlier works, at least the middle of the 5th century. On this point, we agree with A. Szabó ([Szabó1978], p. 65-66) *contra* P. Tannery (*op. cit.*, *supra*, note 15).

[61] As for instance the diagonal and the side of the square.

[62] According to the last line of the passage, an example would be the side of the cube of volume 3: the square built on his side is still incommensurable to the unit (in modern term the side is $\sqrt[3]{3}$ and its square is $\sqrt[3]{9}$ which is incommensurable to 1).

[63] 'καὶ ἀσύμμετροι αἱ μὲν *μήκει* μόνον, αἱ δὲ καὶ *δυνάμει*'. Along with Theodorus' lesson, and contrary to Euclid, the boys did not consider the case of incommensurability in 'power' i.e. lines whose squares are still incommensurable (in modern terms $a$ and $b$ such that the ratio $a^2/b^2$ is still irrational) except briefly and indirectly at the very end of the passage, when Theaetetus speaks of an analogous definition for the side of a cube ('καὶ περὶ τὰ στερεὰ ἄλλο τοιοῦτον', 148b2; in modern terms the cubic roots of integers). Nevertheless, it shows that they were aware of the existence of such lines, and that such magnitudes were well-known at the end of the 5th century.

[64] In modern terms, lines whose lengths are either integers or square roots of integers, with respect to an arbitrary unit of length.

### 4.3. Back to Euclid

Let us return to the definition X.4 of the *Elements*. Two classes are defined. On the one hand, the class of the straight lines commensurable to an assigned straight line, the 'unit', on the other hand the lines incommensurable 'in length' to the same assigned line but commensurable 'in power'. Thus the integers have vanished and the statement is entirely set in geometrical terms. It is a homogeneous geometrical 'definition' and the opposition is the between commensurable and incommensurable magnitudes. On the contrary, according to the boys' definition, in line with Theodorus' lesson, the opposition is between integers and incommensurable lines that entails an amalgamate of two mutually exclusive elements, the discrete and the continuous. Moreover, while the boys definition of 'powers' results from its opposition to 'length', each class in the *Elements* is defined by a common property of its elements,[65] not by exclusion of each other.[66]

Worse still, Theaetetus remarks as in passing that 'powers' are incommensurable to 'length', and in particular to the unit-length. This makes their definition to be mathematically a total failure. Not that the statement is false, it is on the contrary perfectly exact, but it needs a very long demonstration and the proof of many intermediate propositions. In other terms, this alleged definition is definitively not a definition but a deep theorem, using a large part of books VII and VIII of Euclid *Elements*.[67]

**Conclusion**. We have tried to show the importance of the problem of the 'powers' in the second part of Theaetetus' account and its close connection to the whole dialogue, and even to the whole trilogy. However, contrary to the 'Main Standard Interpretation', the fundamental question is not the exact meaning of 'powers', that is 'square' or 'side of square'. Plato's text shows the flaws in the reasoning of Theodorus' brightest pupils.[68] They failed because they considered mainly objects like 'lengths', 'powers' and not relations like 'commensurability/incommensurability in length or in power', though Theodorus' lesson was precisely about the relations of commensurability/incommensurability. The same mistake is at the origin of the failures of all Theaetetus' attempts to define science: from a Protagorean point of view,[69] only sensible objects are worth of interest, not their relations. It is a self-defeating point of view, for the only possibility to study the former are through our senses i.e. through some relations between us and them. It is why, at the end of the dialogue, the '*logos*' (i.e. the 'relation') appears under the guise of two dreams (201c-202c), and the attempts to get its definition ends, once again, in a failure (206c-210b).

---

[65] 'Those magnitudes are said to be commensurable (σύμμετρα') which are measured by the same measure, and those incommensurable ('ἀσύμμετρα') which cannot have any common measure.' (definition 1, Book X).

[66] Definition X.1 considers 'incommensurable' ('*asymmetra*') simply as non-commensurable. However, contrary to Greek and Barbarian, it is not a different name. The prefix 'a' in '*asymmetra*' means indeed that this term is precisely the negation of '*symmetra*', thus both terms are well-defined according to Plato's rules in the *Sophist* and the *Statesman*. However, we certainly do not suggest that Euclids's definitions respected necessarily Plato' criteria of division.

[67] According to the tradition, it will be proved later by Theaetetus (excluding by the way his friend *Socrates*), when he will have become a bright mathematician emancipated from Theodorus' teaching. For details we refer to [Brisson-Ofman2].

[68] As a matter of fact, we have emphasized here the shortcomings in the boys reasoning. But they certainly improved Theodorus' presentation about incommensurability, at least on two fundamental points. First, they consider the general problem, involving infinite quantities, instead of taking 'each case in turn'. And above all, instead of using sensible measures such as the foot, and drawing graphical proofs, they referred to abstract units and abstract demonstrations where figures are simply used as helper for words. This alone would justify Socrates' favorable comment at the end of Theaetetus' account (148b) (cf. also *infra*, note 72).

[69] Protagoras is said to be Theodorus' master ('διδάσκαλος', 179a), and in turn, he is Theaetetus' master; cf. also, *supra*, note 23.

Mathematics is definitively not about sensible objects; it is about measure ('μέτρον') and commensurability, the foundations of any science and technique (*Statesman*, 284a-b). This is certainly not restricted to mathematics. To the objects 'in power' ('δυνάμει') that occupy a central part of the passage, correspond the techniques and sciences as powers, because they are the 'powers' to do something ('to persuade', 'πείθειν', 304c; 'to go to war', 'πολεμητέον', 304e; to judge rightly, 'ὀρθῶς δικάζειν', 305b,c; … ), but also at the end of the *Theaetetus* itself, what Socrates' technique 'has the power to do' (cf. *infra*, end of note 72). Back to the *Theaetetus*, we are able now to understand the huge paradox of Socrates claiming that to *have* the science of the integers is to *have* all of them in one's soul. [70] Faithful to the midwifery, he brings the young boy to clarify the premises put forth by him, [71] which results in the failures of the search of the inquiry. The main obstacle is not, as usually assumed, the complete absence of intelligible forms, but the dreamed hunt of some objects of science, let it be imprints (194c) or birds (197c), while science (as well as technique) is an elusive 'power' of thinking ratios. [72]

## Appendix

For a modern pure logicist, [73] the set of all the integers and the set of the perfect squares (or the set of all even integers) are equivalent for they are all numerically infinite. However, a smarter mathematician would consider the limit of the ratio of the number of their elements less than a given number when this number increases. When this ratio is near to 1 (respectively to 0, respectively is very large), the size of the second set would be more or less the same (respectively much smaller, respectively much bigger) than the first. Let us consider some examples. For any integer $n$

1. There are $n/2$ even integers between 1 and $n$, so that the 'ratio' in the above sense of all integers to the even ones is 2. Thus, there are *twice* 'more' integers than even ones, which is consistent with the intuition.
2. The 'ratio' of all the integers to the perfect squares is greater than the ratio of $n$ over $\sqrt{n}$ (since there are at most $\sqrt{n}$ perfect squares less than $n$) and this ratio is extremely large (since $n/\sqrt{n} = \sqrt{n}$). Thus the second set (of perfect squares) is much 'smaller' than the first one (of all integers).
3. Conversely, if we consider the ratio of all the integers to the 'rectangulars' (i.e. the non-perfect square integers) less than $n$, it is the same as the ratio of $n$ to $(n - \sqrt{n})$, which is very close to 1 when $n$ is very large. [74] Thus there are almost as 'many'

---

[70] Cf. *supra*, §2.1.

[71] 'You forget, my friend, that I (…) claim none of them as mine (…) and am merely acting as a midwife to you' (157c).

[72] The same may be said about Socrates' remark at the end of the passage: 'Excellent my boys' ('Ἄριστά γ' ἀνθρώπων, ὦ παῖδες·', 148b). One has to understand it more or less as: 'Well, now let us see how your master has helped you to think by yourselves', the fundamental purpose of Socrates' maieutic, as asserted at the end of the dialogue. Indeed, the result of the long discussion with Theaetetus is for the latter to 'be pregnant with better thoughts', to 'be less harsh and gentler to your companions, and first and foremost, to 'have the wisdom not to think you know that which you do not know' ('σωφρόνως οὐκ οἰόμενος εἰδέναι ἃ μὴ οἶσθα', 210c); then, Socrates to conclude: 'This is all my art *has the power to do*' ('τοσοῦτον γὰρ μόνον ἡ ἐμὴ τέχνη δύναται'). See however also, *supra*, note 68.

[73] By this term, we do not mean the modern logician, but someone, more often than not a philosopher, considering any work, ancient or modern, from the point of view of (elementary) logics.

[74] Since $n/(n - \sqrt{n}) = 1/(1 - (\sqrt{n}/n))$ and $\sqrt{n}/n = 1/\sqrt{n}$ is very small when $n$ is very large. For instance, for $n = 10000$, we get: $n/(n - \sqrt{n}) = 10000/9900 \approx 1,01$.

'rectangular' integers as integers, that is consistent with both the result in point 2. above and, once again, the intuition.

Though it is a modern presentation, this is probably what Plato, and the mathematicians of the 5th-4th centuries had in mind, since in ancient Greek, one of the meanings of 'ἄπειρον' was something 'very large' (cf. for example the meaning of *half* the 'whole integer' ('τοῦ ἀριθμοῦ ἅπας'), in *Phaedo*, 104a-b). Indeed, the modern logicist will distinguish between finite and enumerable infinite, and all the numerable infinites are the same for him. For the smarter mathematician, there is a measure for the different enumerable infinites. If our analysis is correct, Plato gives here, inside a mathematical background and in the extreme case of infinite classes, a first introduction on the two kinds of comparison detailed in the *Statesman*. One considers on the one hand, in agreement with the logicist, the pure relation of the more and less (283c-d), on the other hand, with the smarter mathematician, the existence of a measure. The latter opens the way to the possibility to decide what is correct, and what is not, in different situations i.e. to the 'just measure' (283e-284c). Anyway, even in Plato's time, such considerations were only within the scope of someone well-skilled in mathematics. This may explain the difficulties for the understanding of this passage even in the (late) Antiquity [75] rather than some alleged mathematical weakness of his author [76] contradicted by all Plato's writings.

# Bibliography

[Ackrill1997]: John Ackrill, *Essays on Plato and Aristotle*, Oxford Univ. Press, p. 93-109, 1997

[Allman1877]: George Allman, *Greek Geometry from Thales to Euclid*, Dublin Univ. Press, 1877

[Bernadete2007]: Seth Bernadete, *The Being of the Beautiful: Plato's* Theaetetus*, Sophist*, and Statesman, Chicago Univ. Press, 2007

[Brisson2003]: Luc Brisson (ed.), *Platon-Le politique*, Flammarion, 2003

[Brisson-Ofman1]: Luc Brisson-Salomon Ofman, "Theodorus' lesson in Plato's *Theaetetus* (147d1-d6) Revisited-A New Perspective", to appear
[Brisson-Ofman2]: Luc Brisson-Salomon Ofman, "The Philosophical Interpretation of Plato's *Theaetetus* and the Final Part of the Mathematical Lesson (147d7-148b)", to appear

[Brown2010]: Lesley Brown, "Definition and Division in Plato's *Sophist*", in *Definition in Greek Philosophy*, D. Charles ed., Oxford Univ. Press, 2010 , p. 151-171

[Burnett1899]: John Burnett (ed.), *Platonis Opera,* t. 1, Clarendon Press, 1899

---

[75] We refer for instance to the *Anonymous commentary* of the *Theaetetus* (cf. [Diels1905]).
[76] Cf. *supra*, note 15.

[Burnyeat 1978]: Myles Burnyeat, "The Philosophical Sense of Theaetetus' Mathematics", *Isis*, 69, 1978, 489-514

[Caveing1998]: Maurice Caveing, *L'Irrationalité dans les mathématiques grecques jusqu'à Euclide*, vol. 3, Septentrion, 1998

[Diels1905]: Herman Diels, *Anonymer Kommentar zu Platons* Theaetet, Papyrus 9782. Nebst drei Bruchstücken philosophischen Inhalts, Pap. N. 8; P. 9766. 9569, unter Mitwirkung von J.L. Heiberg, bearbeitet von H. Diels und W. Schubart, Berlin, 1905 (new edition by, Guido Bastiani and David Sedley, "Commentarium in Platonis *Theaetetum*". in *Corpus dei papyri filosofici greci e latini*, edition, translation, notes, Part III, 1995, p. 227-562)

[Fowler1921]: Harold Fowler (ed.), *Plato in Twelve Volumes* (*Theaetetus, Sophist, Statesman*), vol. 12, Harvard University Press, 1921

[Fowler1999]: David Fowler, *The Mathematics of Plato's Academy: A New Reconstruction*, Clarendon Press, 1999

[Gill2010]: Mary Louise Gill, "Division and Definition in Plato's *Sophist* and *Statesman*", in *Definition in Greek Philosophy*, D. Charles ed., Oxford Univ. Press, 2010, p. 172-199
[Gill2016]: Mary Louise Gill, "Method and Metaphysics in Plato's *Sophist* and *Statesman*", *The Stanford Encyclopedia of Philosophy*, Edward N. Zalta (ed.), Winter 2016 Edition

[Guyer-Wood1998]: Paul Guyer-Allen Wood (ed.), *Kant-Critique of pure reason*, Cambridge Univ. Press, 1998,

[Harte2002]: Verity Harte, Plato, *On parts and wholes*, Oxford Univ. Press, 2002

[Kahane1985]: Jean-Pierre Kahane, "la théorie de Théodore des corps quadratiques réels", *L'enseignement mathématique*, 31, 1985, p. 85-92

[Kahn2007]: Charles Kahn: "Why Is the Sophist a Sequel to the Theaetetus?", *Phronesis,* 52 (2007), p. 33-57

[Klein1977]: Jacob Klein (ed.), Plato's Trilogy: "Theaetetus", "The Sophist" and "The Statesman", Chicago Univ. Press, 1977

[Knorr1975]: Wilbur Knorr, *The evolution of the Euclidean elements*, Reidel, 1975

[Lawlord1979], Robert and Deborah Lawlor, *Mathematics useful for understanding Plato by Theon of Smyrna Platonic Philosopher*, Wizards Bookshelf, 1979

[Ofman2010]: Salomon Ofman, "Une nouvelle démonstration de l'irrationalité de racine carrée de 2 d'après les Analytiques d'Aristote", *Philosophie antique*, 10, 2010, p. 81-138

[Ofman2014]: Salomon Ofman, "Comprendre les mathématiques pour comprendre Platon-*Théétète* (147d-148b)", *Lato Sensu*, I, 2014, p. 70-80

[Ofman2018]: Salomon Ofman, "Irrationalité, désordre, géométrie", in *Complexité et Désordre-Rencontres*, Matériologique, 2018, to appear

[Rowe]: Christopher Rowe (ed.), *Plato, Theaetetus and Sophist*, Cambridge Univ. Press, 2015

[Ryle1939]: Gilbert Ryle, "Plato's Parmenides (ii)", *Mind*, 48, 1939, p. 302–325

[Sayre2006]: Kenneth Sayre, *Metaphysics and Method in Plato's Statesman*, Cambridge Univ. Press, 2006

[Szabó1978]: Arpád Szabó, *The Beginning of Greek mathematics*, transl. Umgar, Reidel, 1978

[Tannery1876]: Paul Tannery, *Mémoires scientifiques*, 17 vol., J. Heiberg and H. Zeuthen (eds), Gauthier-Villars, 1912-1950

[Unguru2002]: Sabetai Unguru, "Amicus Plato sed…: Fowler's New Mathematical Reconstruction of the Mathematics of Plato's Academy", *Annals of Science*, 59, 2002, p. 201–210